\documentclass[11pt]{tran-l}
 \usepackage{array}
 \usepackage{tikz}
 \usetikzlibrary{matrix,arrows}

\usepackage[centertags]{amsmath}
\usepackage{amsfonts}
\usepackage{amssymb}
\usepackage{amsthm}
\usepackage{graphics, graphicx}
\usepackage{newlfont}
\usepackage[all]{xy}
\usepackage{url}
\theoremstyle{plain}
\newtheorem{thm}{Theorem}[section]
\newtheorem{cor}[thm]{Corollary}
\newtheorem{lem}[thm]{Lemma}
\newtheorem{prop}[thm]{Proposition}
\theoremstyle{definition}
\newtheorem{defn}{Definition}[section]
\theoremstyle{remark}

\theoremstyle{Example}
\newtheorem{exmp}{Example}[section]
\theoremstyle{Conjecture}

\numberwithin{equation}{section}

\begin{document}

\title[Pascal Determinantal Arrays and Log-Concavity]
{Khayyam-Pascal Determinantal Arrays, Star of David Rule and Log-Concavity}
\author{Hossein Teimoori Faal and Hasan Khodakarami}
\address
{Department of
Computer Science, Faculty of Mathematical and Computer Sciences, 
Allameh Tabataba'i University, Tehran, Iran.	
}
\email{hossein.teimoori@atu.ac.ir}
\address
{Ministry of Information and Communication Technology (ICT), Zanjan, Iran.}
\email{hasan.khodakarami@gmail.com}

\maketitle

%%% ----------------------------------------------------------------------

\begin{abstract}
	
In this paper we we develop a new geometric method to answer the log-concavity questions related to a nice class of combinatorial sequences arising from the Pascal triangle.

\end{abstract}

%%% ----------------------------------------------------------------------

\section{Introduction}

One of the important task in \textit{enumerative combinatorics} is to determine \textit{log-concavity} of a combinatorial sequence.
\begin{defn}
A sequence $a_{0}, a_{1},\ldots , a_{n}$ of real numbers is said to be \textit{concave} if
$\frac{a_{i-1}+a_{i+1}}{2}\leq a_{i}$ for all $1 \leq i \leq n-1$, and \textit{logarithmically concave}
(or log-concave for short) if $a_{i-1}a_{i+1}\leq a_{i}^{2}$ for all $1 \leq i \leq n-1$.
\end{defn}
\begin{defn}
The sequence $a_{0}, a_{1},\ldots , a_{n}$ is called \textit{symmetric} if $a_{i}=a_{n-i}$ for $0 \leq i \leq n$.
\end{defn}
\begin{defn}
We say that a polynomial $a_{0} + a_{1} q + \cdots + a_{n} q^{n}$ has a certain property (such as log-concave or symmetric) if its sequence $a_{0}, a_{1},\ldots , a_{n}$ of coefficients has the property.
\end{defn}
There are many ways to prove the log-concavity of a combinatorial sequence. One of the classic method of proof is a direct combinatorial approach, which is of significant interest for combinatorial people.
\begin{exmp}
The best-known log-concave sequence is the $n$-th row of
\textit{Khayyam-Pascal's triangle}:
\\
\begin{center}
${ {n\choose 0} , {n\choose 1} , {n\choose 2} , \ldots , {n\choose n} }$.
\end{center}

Here, the log-concavity is easy to show directly because of the explicit formula ${n \choose k} = \frac{n!}{k!(n-k)!}$. Indeed,
$$
\frac{{n \choose k}^{2} }{{n \choose k-1}{n \choose k+1}}=\frac{(k+1)(n-k+1)}{k(n-k)} > 1,
$$
which is equivalent to $n>-1$ (or $n\geq 0$), as required.
\end{exmp}
\begin{exmp}
For the sequence of the $n$-th diagonal of
the Khayyam-Pascal triangle:
\begin{center}
${ {n\choose 0} , {n+1\choose 1} , {n+2\choose 2} , \ldots , {n+k\choose k} } , \ldots$,
\end{center}
again, we have
\begin{center}
$\frac{{n+i \choose i}^{2} }{{n+i-1 \choose i-1}{n+i+1 \choose i+1}}=\frac{(n+i)(i+1)}{i(n+i+1)} > 1$,
\end{center}
which is equivalent to $n>0$.
\end{exmp}

In spite of the \textit{geometric idea} behind the definition of the log-concavity of a sequence, to the best of our knowledge, there is no geometric approach to tackle this issue. In this paper, we develop a new geometric method to answer the log-concavity questions related to a nice class of combinatorial sequences arising from the Khayyam-Pascal triangle.

\section{Khayyam-Pascal Array and Parallelepiped Determinantal Identities}

\begin{figure}[h]
\begin{center}
\includegraphics[scale=.7]{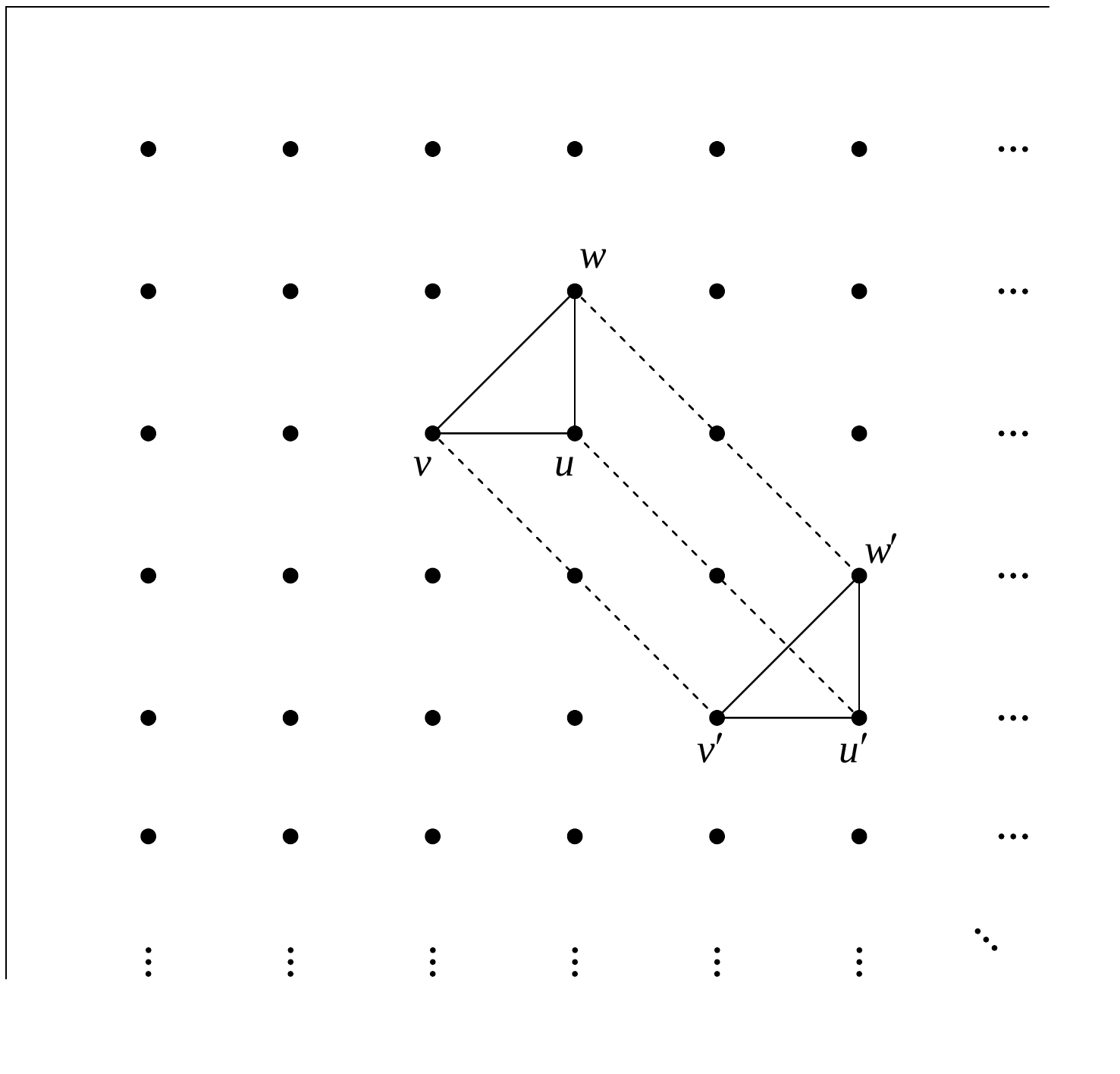}\\
\end{center}
\caption{
A Determinantal Parallelepiped
\label{fig:F1}
}
\end{figure}

Consider a $45\,^{\circ}$ rotation of the Khayyam-Pascal triangle which we call it
\textit{Khayyam-Pascal squared array} \cite{m1}. Now, we construct a parallelepiped with two triangles as its bases which is shown with six entries of
this array and the corresponding edges in 
Figure~\ref{fig:F1}.

Then, we have the following determinantal identities which are the direct consequence of the recurrence relation
for the Khayyam-Pascal array. 

\begin{prop}(Parallelepiped Determinantal Identities)
\begin{eqnarray}
i) \quad \left|
\begin{array}[c]{ccc}
u&v \\
u'&v'\\
\end{array}
\right|
&=&\,
\left|
\begin{array}[c]{ccc}
w&v \\
w'&v'\\
\end{array}
\right|,\nonumber\\
ii) \;\; \left|
\begin{array}[c]{ccc}
w&v \\
w'&v'\\
\end{array}
\right|
&=&
\left|
\begin{array}[c]{ccc}
w&u \\
w'&u'\\
\end{array}
\right|.\nonumber
\end{eqnarray}
In other words, the determinants formed by three faces of the parallelepiped $uvwu'v'w'$ in 
Figure~\ref{fig:F1} are equal.
\end{prop}
\begin{proof}
By the rule of Khayyam-Pascal array, we have
\begin{center}
$u=v+w$
\\
$u'=v'+w'$.
\end{center}
Now, multiplying the above equalities by $v$ and $v'$, respectively, we get
\\
\begin{center}
$uv'=vv'+wv'$
\\
$u'v=vv'+w'v$.
\end{center}
Subtracting the above equalities, we obtain
\begin{center}
$uv' - u'v = wv' - w'v$,
\end{center}
or equivalently
\\
\begin{center}
$
\left|
\begin{array}[c]{ccc}
u&v \\
u'&v'\\
\end{array}
\right|
=
\left|
\begin{array}[c]{ccc}
w&v \\
w'&v'\\
\end{array}
\right|,
$
\end{center}
which is the first determinantal identity. The second one can be proved in a similar way
and left to the reader as a simple exercise.
\end{proof}
\begin{prop}
Every diagonal of the Khayyam-Pascal triangle is log-concave.
\end{prop}
\begin{proof}

First of all note that the diagonals of the Khayyam-Pascal triangle correspond to the columns (rows) of the Khayyam-Pascal squared array. Now, we use the previous determinantal identities in their special cases to give a new geometric proof of the log-concavity of the diagonals of the Khayyam-Pascal triangle. 

\newpage 

To this end, consider three consecutive terms
$a_{k-1}, a_{k}, a_{k+1}$ in any arbitrary column of the Khayyam-Pascal squared array, as shown in Figure~\ref{fig:F2}.

\begin{figure}[h]
	\begin{center}
		\includegraphics[scale=.8]{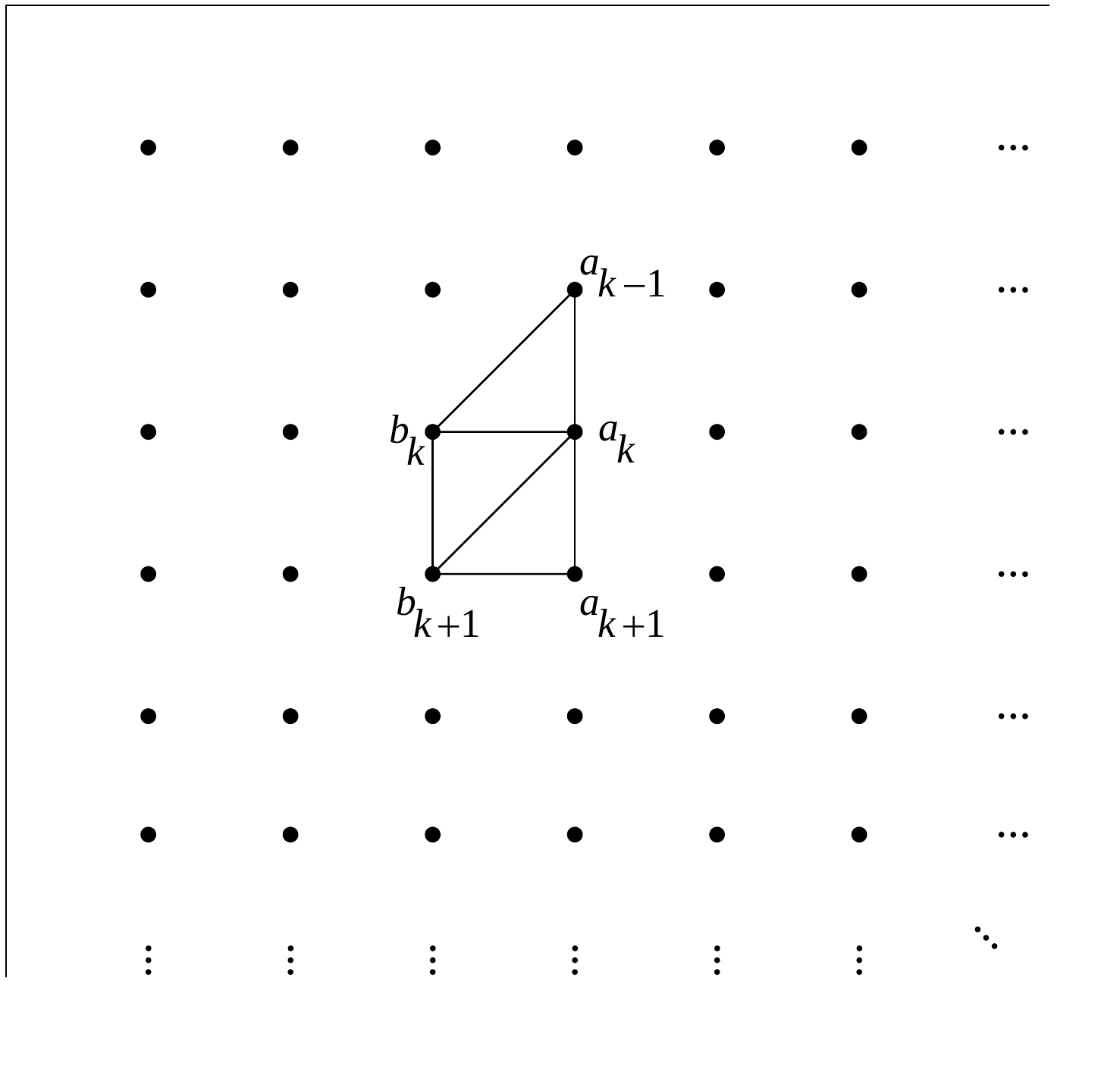}
	\end{center}
	\caption{
		Log-Concavity of Diagonals of the Khayyam-Pascal Triangle
		Array.
		\label{fig:F2}
	}
\end{figure}

We consider a parallelepiped in its special case where two antipodal vertices ($u$ and $w'$ in Figure~\ref{fig:F1}) coincide.
Here, those vertices correspond to two equal entries $a_{k}$. By Proposition $2.1$, we have

\begin{center}
$
\left|
\begin{array}[c]{ccc}
a_{k}&a_{k+1} \\
a_{k-1}&a_{k}\\
\end{array}
\right|
=
\left|
\begin{array}[c]{ccc}
b_{k}&a_{k} \\
b_{k+1}&a_{k+1}\\
\end{array}
\right|.
$
\end{center}
But, we already know that the 2-by-2 determinant in the right-hand side of the above identity is a Narayana number $[4]$. Therefore, we obtain
\\
\begin{center}
$
\left|
\begin{array}[c]{ccc}
a_{k}&a_{k+1} \\
a_{k-1}&a_{k}\\
\end{array}
\right|
\geq 0,
$
\end{center}
and this completes the proof.
\end{proof}

Next we prove the log-concavity of the rows of the Khayyam-Pascal triangle, using the same technique.
\begin{prop}
Every row of the Khayyam-Pascal triangle is log-concave.
\end{prop}
\begin{proof}
We note that the rows of the Khayyam-Pascal triangle correspond to the diagonals of the the Khayyam-Pascal squared array. Consider an special parallelepiped $vuwv'u'v$, as shown in Figure~\ref{fig:F3}.

\begin{figure}[h]
\begin{center}
\includegraphics[scale=.8]{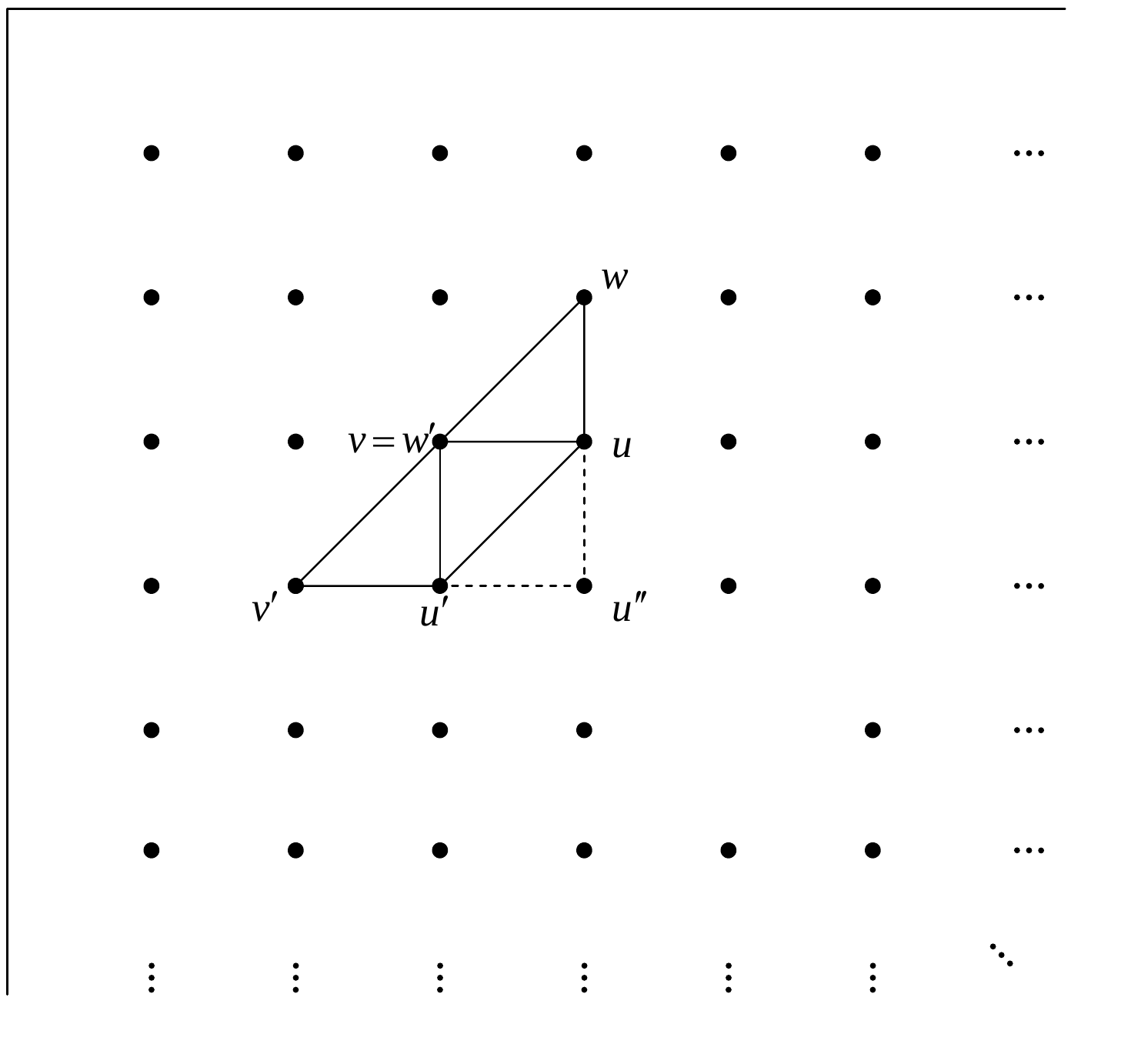}
\end{center}
\caption{
Log-Concavity of Rows of Khayyam-Pascal Triangle
Array.
\label{fig:F3}
}
\end{figure}

Then, we have

\begin{center}
$
\left|
\begin{array}[c]{ccc}
v&w \\
v'&v\\
\end{array}
\right|
=
\left|
\begin{array}[c]{ccc}
w'&w\\
u'&u\\
\end{array}
\right|.
$
\end{center}
On the other hand, from the parallelepiped $w'uwu'u''u$ we get
\begin{center}
$
\left|
\begin{array}[c]{ccc}
w'&w \\
u'&u\\
\end{array}
\right|
=
\left|
\begin{array}[c]{ccc}
w'&u\\
u'&u''\\
\end{array}
\right|.
$
\end{center}
Therefore, we conclude that

\begin{center}
$ v^{2}-wv'
= \left|
\begin{array}[c]{ccc}
w'&u\\
u'&u''\\
\end{array}
\right|. $
\end{center}                  
But, again the last determinant in the above equality is the Narayana number and a non-negative integer.
This completes the proof.
\end{proof}
\begin{defn}
We call an array a \textit{row log-concave} (diagonal log-concave) array, if every row (diagonal) of this
array is log-concave.
\end{defn}
As in the paper of McNamara and Sagan \cite{m2}  for every array $A=(a_{ij})_{i,j\geq 0}$, we will call the determinants
$
\left|
\begin{array}[c]{ccc}
a_{i,j}&a_{i,j+1}\\
a_{i+1,j}&a_{i+1,j+1}\\
\end{array}
\right|
$,
its \textit{adjacent minors}. From the proofs of the two previous propositions, we get the following
interesting result.
\begin{cor}
Every diagonal log-concave array with non-negative adjacent minors, is also a row log-concave array.
\end{cor}

\section{Khayyam-Pascal Determinantal Arrays}

In this section, we introduce an infinite class of arrays of numbers as a generalization of the standard
Khayyam-Pascal squared array. We will denote the entries of the the Khayyam-Pascal squared array by
$
P=
\big( 
P_{i,j} =
{i+j\choose i}
\big)_
{i,j
\geq 0}
$.
Our main goal here is to prove that
the members of this new class of arrays are diagonal and row log-concave, again using geometric ideas.

\begin{figure} 
\begin{center}
\includegraphics[scale=1.2]{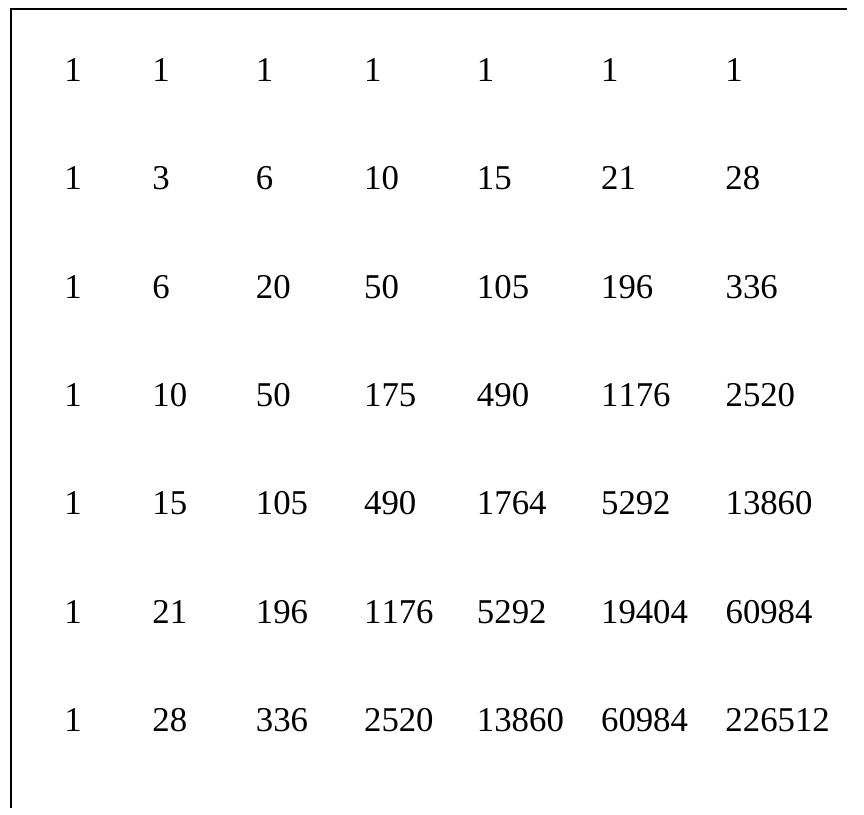}
\end{center}
\caption{
Khayyam-Pascal Determinantal array of order 2.
}
\label{fig:F4}
\end{figure}

\begin{defn}
A Khayyam-Pascal determinantal array of order $k$, $k\geq 1$, is an infinite array 
$
PD_{k} =
\big(
P_{i,j}^{(k)}
\big)_{i,j \geq 0}
$
in which $P_{i,j}^{(k)}$ 
is the determinant of a $k$-by-$k$ subarray of
the Khayyam-Pascal squared array, starting form $(i,j)$-entry. Namely,
\begin{center}
$ P_{i,j}^{(k)}:= \left|
\begin{array}[c]{ccc}
P_{i,j}&\ldots & P_{i,j+k-1} \\
\vdots&\ddots & \vdots \\
P_{i+k-1,j}&\ldots & P_{i+k-1,j+k-1}
\end{array}
\right|.
$
\end{center}
\end{defn}
\begin{exmp}
A Khayyam-Pascal determinantal array of order 2 has shown in Figure~\ref{fig:F4}. This is a well-known array
which is the squared-form of the so-called \textit{Narayana triangular array} (see A001263 in \cite{m3}).
\end{exmp}

\newpage 

\begin{figure}[h]
\begin{center}
\includegraphics[scale=.80]{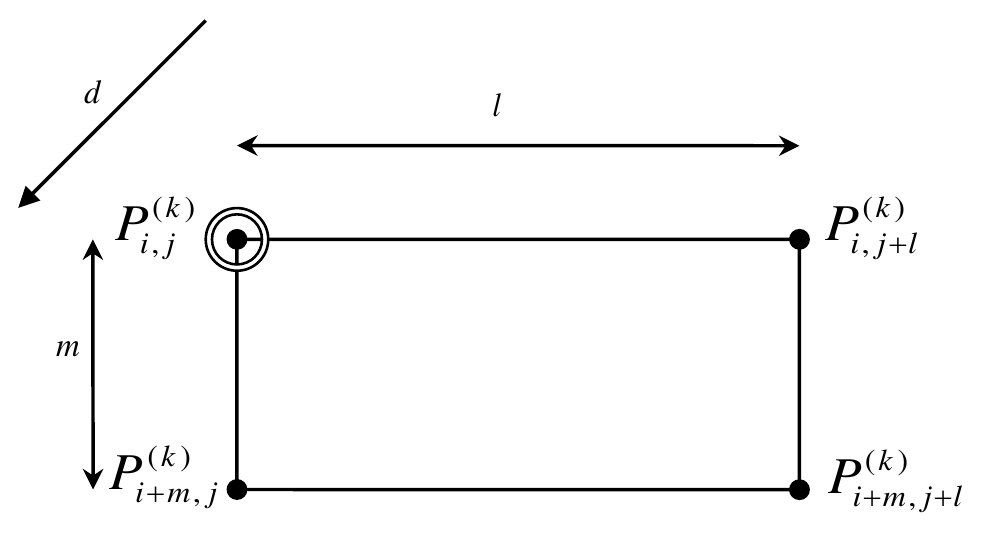}
\end{center}
\caption{
 Weighted Version of Star of David.
\label{fig:F5}
}
\end{figure}

In \cite{m4}, the authors have shown that if we define the weight of any arbitrary rectangle whose vertices are
the entries of the Khayyam-Pascal determinantal array of order $k$ as shown in Figure~\ref{fig:F5}, by

\begin{center}

$W:=\frac{P_{i+m,j+l}^{(k)}.
P_{i,j}^{(k)}}
{P_{i+m,j}^{(k)}.P_{i,j+l}^{(k)}},$

\end{center}

%\begin{equation*}\label{StarDavidRule1}
%\xymatrix{
%&\ar[dl]_(.5){d}       &             &            %&              \\
%&*+[o][F-]{\bullet}\ar @{-}[rrr]^(.5){l}^<<{P_{i,j}^{(k)}}^>{P_{i,j+l}^{(k)}}   &  &    &{\bullet}\ar@{-}[dd]              \\
%&              &             &            &                      \\
%&{\bullet}\ar @{-}[rrr]_<{P_{i+m,j}^{(k)}}_>{P_{i+m,j+l}^{(k)}} \ar@{-}[uu]^(.5){m}  &     &     & {\bullet}             }
%\end{equation*}

then when we move the anchor, the circled-vertex, along the diagonal of the Khayyam-Pascal determinantal
array (indicated by the arrow $d$ in Figure~\ref{fig:F5}), the weights remain unchanged. They called this property
the \textit{weighted-version of the Star of David Rule}. As they have shown in another paper \cite{m5}, the
weighted-version of the Star of David Rule can also be used to prove the following interesting property of this
new class of arrays.

\begin{figure} [h]
\begin{center}
\includegraphics[scale=.8]{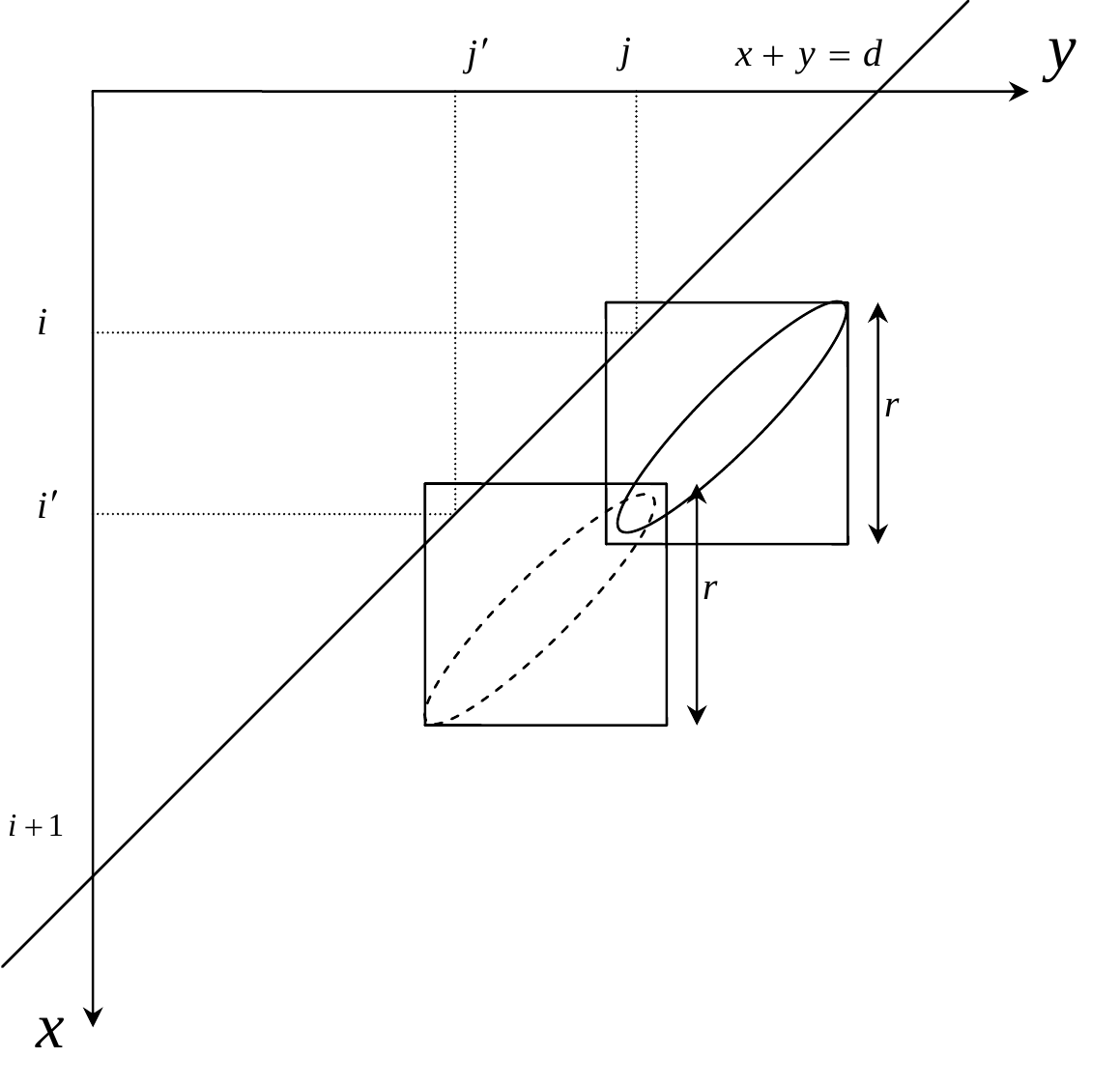}
\end{center}
\caption{ 
Ratio of Determinants in Khayyam-Pascal Determinantal array
}
\label{fig:F6}
\end{figure}

\begin{prop}
In any Khayyam-Pascal determinantal array, the ratio of any pair of $r$-by-$r$ minors along any arbitrary diagonal
$x+y=d$ of the array is the same as the ratio of the product of the entries appearing in their back diagonals
parallel to $d$ (see Figure \ref{fig:F6}). In other words, we have
\begin{center}
$\frac
{\left|
\begin{array}[c]{ccc}
P_{i,j}^{(k)}&\ldots & P_{i,j+r-1}^{(k)} \\
\vdots&\ddots & \vdots \\
P_{i+r-1,j}^{(k)}&\ldots & P_{i+r-1,j+r-1}^{(k)}
\end{array}
\right|}{\left|
\begin{array}[c]{ccc}
P_{i', j'}^{(k)}&\ldots & P_{i',j'+r-1}^{(k)} \\
\vdots&\ddots & \vdots \\
P_{i'+r-1,j'}^{(k)}&\ldots & P_{i'+r-1,j'+r-1}^{(k)}
\end{array}
\right|}
=
\frac
{P_{i,j+r-1}^{(k)} \ldots P_{i,j+r-1}^{(k)} }{P_{i',j'+r-1}^{(k)} \ldots P_{i',j'+r-1}^{(k)} }
$
\end{center}
\end{prop}

\newpage	
	
The following lemma is the key in the proof of diagonal log-concavity of the Khayyam-Pascal determinantal Arrays.
\begin{lem}
For every integer $n\geq 1$, the log-concave sequence $\{a_{i}\}_{i\geq 1}$
satisfies the following inequality
$$
\frac{a_{2}a_{n+1}}{a_{1}a_{n+2}} \geq 1.
$$
\end{lem}
\begin{proof}
We use induction on $n$. The basis case, $n=1$, is just the definition of the log-concavity of the sequence $\{a_{i}\}_{i \geq 1}$. Now, let us assume by induction hypothesis that the assertion is true for $n-1$.
Hence, we have
\begin{center}
$$
1 \leq \frac{a_{2}a_{n}}{a_{1}a_{n+1}} = (\frac{a_{2}a_{n}}{a_{1}a_{n+1}})(\frac{a_{n+1}a_{n+2}}{a_{n+1}a_{n+2}}) = (\frac{a_{2}a_{n+1}}{a_{1}a_{n+2}})(\frac{a_{n}a_{n+2}}{a_{n+1}^{2}}).
$$
\end{center}
Thus, we get
$$
\frac{a_{2}a_{n+1}}{a_{1}a_{n+1}} \geq \frac{a_{n+1}^{2}}{a_{n}a_{n+2}} \geq 1.
$$
The later inequality holds because of the definition of the log concavity of the sequence $\{a_{i}\}_{i \geq 1}$.
This completes the proof by induction.
\end{proof}
Now, we are at the position to state our main result of this section.
\begin{thm}
For every integer $k\geq 1$, the Khayyam-Pascal determinantal array of order $k$ is diagonal log-concave.
\end{thm}
\begin{proof}
Assume that $\alpha$, $\beta$, $\theta$, $\gamma$ are four entries of the
Khayyam-Pascal determinantal array of order $k$ such that $\beta$, $\theta$, $\gamma$ are three consecutive
diagonal entries, as shown in Figure~\ref{fig:F7}.

\begin{figure} [h]
\begin{center}
\includegraphics[scale=.6]{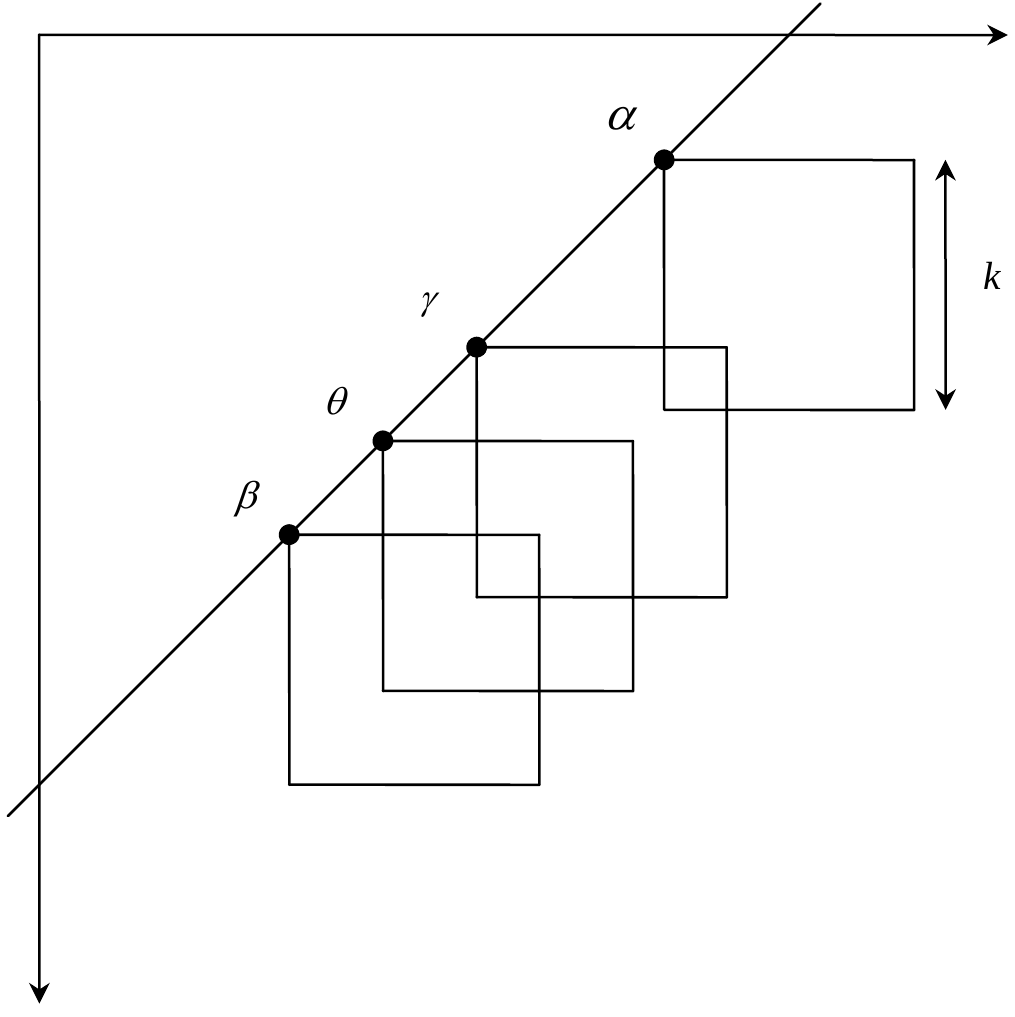}
\end{center}
\caption{ 
Four Entries of A diagonal of Khayyam-Pascal Determinantal Array
}
\label{fig:F7}
\end{figure}

\newpage

Clearly the back diagonal entries of these four entries of the Khayyam-Pascal determinantal array of order $k$, as the four $k$-by-$k$ minors of the Khayyam-Pascal squared array, lie in some diagonal of the Khayyam-Pascal squared array. For simplicity of arguments, we will show their entries from south-west to north-east by ${\beta}_{1}$, ${\beta}_{2}$, $\ldots$, ${\beta}_{k}$,  ${\theta}_{1}$, ${\theta}_{2}$, $\ldots$, ${\theta}_{k}$,  ${\gamma}_{1}$, ${\gamma}_{2}$, $\ldots$, ${\gamma}_{k}$ and ${\alpha}_{1}$, ${\alpha}_{2}$, $\ldots$, ${\alpha}_{k}$, respectively. It is not hard to see that we have the following relations among their
entries:
\begin{eqnarray}
&&\beta_{2}=\theta_{1}  , \beta_{3}=\theta_{2}, \ldots, \beta_{k}=\theta_{k-1},\nonumber\\
&&\theta_{1}=\gamma_{1} , \theta_{3}=\gamma_{2},\ldots, \theta_{k}=\gamma_{k-1}.\nonumber
\end{eqnarray}
To prove the log-concavity, it suffices to show that ${\theta}^{2} - \beta \gamma \geq 0$. But,
using the determinants ratio Proposition $3.1$ and the above relations, we have
\begin{eqnarray}
&&{\theta}^{2} - \beta \gamma \nonumber\\
&&= \left(\frac{\alpha}{\alpha_{1}\cdots\alpha_{k-1}\alpha_{k}}\right)^{2}
\left[(\theta_{1}\cdots\theta_{k-1}\theta_{k})^{2}-(\beta_{1}\cdots\beta_{k-1}\beta_{k})
(\gamma_{1}\cdots\gamma_{k-1}\gamma_{k})\right],\nonumber\\
&&= \left(\frac{\alpha}{\alpha_{1}\cdots\alpha_{k-1}\alpha_{k}}\right)^{2}
\left[(\beta_{2}\beta_{3}^{2}\cdots\beta_{k}^{2}\gamma_{k-1})(\beta_{2}\gamma_{k-1}-\beta_{1}\gamma_{k}        \right].\nonumber
\end{eqnarray}
Therefore we need to prove that $\frac{\beta_{2}\gamma_{k-1}}{\beta_{1}\gamma_{k}}\geq 1$, which
is nothing more than the inequality of the key lemma, Lemma $3.2$, by the \emph{renaming technique}.
\end{proof}
Next, we prove the row log-concavity of the Khayyam-Pascal determinantal array.
\begin{thm}
For every integer $k\geq1$, the Khayyam Pascal determinantal array of order $k$ is a row log-concave array.
\end{thm}

\begin{proof}
Using Corollary $2.4$, it is only suffices to prove that every adjacent minor of the Khayyam-Pascal determinantal array of order $k$ is nonnegative. Now by the Proposition $3.1$ about the ratio of determinants along
the diagonal $x+y=d$, we get
$$
\frac{\left|
\begin{array}[c]{ccc}
P_{i,j}^{(k)}&P_{i,j+1}^{(k)} \\
P_{i+1,j}^{(k)}&P_{i+1,j+1}^{(k)}\\
\end{array}
\right|}
{\left|
\begin{array}[c]{ccc}
1&P_{i+j,1}^{(k)} \\
1 &P_{i+j+1,1}^{(k)}\\
\end{array}
\right|}=
\frac{P_{i+1,j}^{(k)}P_{i,j+1}^{(k)}}{P_{i+j,1}^{(k)}},
$$
which is clearly a positive integer. Thus, to prove that the adjacent minor
$
\left|
\begin{array}[c]{ccc}
P_{i,j}^{(k)}&P_{i,j+1}^{(k)} \\
P_{i+1,j}^{(k)}&P_{i+1,j+1}^{(k)}\\
\end{array}
\right|$ is a nonnegative integer,
we only need to show that
$\left|
\begin{array}[c]{ccc}
1&P_{i+j,1}^{(k)} \\
1 &P_{i+j+1,1}^{(k)}\\
\end{array}
\right|$
is positive for every $i,j\geq 0$, which is equivalent to show that the first column, starting form $0$, of the
Khayyam-Pascal determinantal array of order $k$ is an increasing sequence. It is not hard to see that
this first column is indeed the $k$th column of the Khayyam-Pascal squared array \cite{m1}. Finally we need
to show that for every $l\geq 0$, we have
$$
\frac{{l+k\choose k}}{{(l-1)+k\choose k}}>1,
$$
which is equivalent to inequality $k>0$ or $k\geq 1$, as required.
\end{proof}

\end{document}